\documentclass[12pt]{amsart}
\usepackage{amssymb,amsmath,amsthm,mathrsfs,multirow,enumerate}

\usepackage{graphicx} 
\usepackage[all]{xy}
\usepackage[bottom=1in, top=1in, left=1in, right=1in]{geometry}

\numberwithin{equation}{section} 

\theoremstyle{plain}
\newtheorem{thm}{Theorem}
\newtheorem{cor}[thm]{Corollary}
\newtheorem{lemma}[thm]{Lemma}

\theoremstyle{definition}

\newcommand{\bbm}{\begin{bmatrix}}
\newcommand{\ebm}{\end{bmatrix}}
\newcommand\norm[1]{\left\lVert#1\right\rVert}

\begin{document}
\title[The triangle inequality for graded real vector spaces] {The triangle inequality for graded real vector spaces} 
\author{Songpon Sriwongsa and Keng Wiboonton}

\address{Songpon Sriwongsa \\ (1) Department of Mathematics \\ Faculty of Science \\ King Mongkut's University of Technology Thonburi (KMUTT) \\ Bangkok 10140, Thailand \\ \newline (2) Center of Excellence in Theoretical and Computational Science (TaCS-COE) \\ Science Laboratory Building \\ Faculty of Science \\ King Mongkut's University of Technology Thonburi (KMUTT) \\ Bangkok 10140, Thailand}
\email{\tt songpon.sri@kmutt.ac.th, songponsriwongsa@gmail.com}

\address{Keng Wiboonton\\ Department of Mathematics and Computer Science \\ Faculty of Science\\ Chulalongkorn University \\Bangkok 10330, Thailand}
\email{\tt kwiboonton@gmail.com}

\begin{abstract}
	In this paper, we prove that a natural candidate for a homogeneous norm on a graded Lie algebra of any length satisfies the triangle inequality which answers Moskowitz's question in \cite{Moskowitz2}. 
\end{abstract}
\maketitle

\bigskip

\section{Introduction}
In \cite{Moskowitz1}, M. Moskowitz extended the classical theorem of Minkowski on lattice points and convex bodies in $\mathbb{R}^n$ to simply connected nilpotent Lie groups $G$ with a  $\mathbb{Q}$-structure whose Lie algebra $\mathfrak{g}$ admits a length 2 grading. For this purpose, the author proved the triangle inequality for a certain natural homogeneous norm of the Lie algebra associated with the grading of length 2. Later, in \cite{Moskowitz2}, the same author extended the homogeneous norms for which the triangle inequality holds to gradings of length 3 and 4. At the end of the paper \cite{Moskowitz2}, the author demonstrated that the method of his proof cannot be carried on to the case of length 5 and higher, so it was left as an open problem. However, the arguments in the proof contains some errors. In \cite{Moskowitz2}, the inequality (13) does not imply the inequality (12) in general. In fact, the inequalities (6) and (12) do not hold for length $r = 2, 3,$ or $4$. One can easily find a counterexample for (12) by taking $r = 4, \norm{X'} = 2$ and $\norm{Y'} = 1$.

In this paper, we claim that for the grading of any length $r \ge 1$, the homogeneous norms defined as in \cite{Moskowitz2} satisfy the triangle inequality. Our proof is quite elementary. The main tools in our proof are only H\"older's inequality and Muirhead's inequality. We refer the reader to \cite{CKP, Steele} for the latter inequality.

A Lie algebra $\mathfrak{g}$ is said to be \textit{graded} if there is a finite family of subspaces $V_1, \dots, V_r$ with $\mathfrak{g} \ = \ V_1 \oplus \dots \oplus V_r$ satisfying $[V_i, V_j] \subseteq V_{i+j}$ for all $i, j$ (here, $V_k = 0$ if $k > r$). The integer $r$ is its \textit{length}. For a graded Lie algebra $\mathfrak{g}$, we define, for $t \in \mathbb{R}^\times$,
$$\alpha_t(v_1, \dots, v_r) \ = \ (tv_1, t^2v_2, \dots, t^rv_r).$$ 
Then $\alpha_t$ is Lie algebra automorphism of $\mathfrak{g}$. Moreover, if $\mathfrak{g}$ is a graded Lie algebra, then $\mathfrak{g}$ must be nilpotent (see \cite{Moskowitz1} and \cite{Moskowitz2}). 

As in \cite{Moskowitz1} and \cite{Moskowitz2}, we define a \textit{homogeneous norm} on a graded Lie algebra $\mathfrak{g}$ as a function $ \norm{\cdot} :\mathfrak{g} \rightarrow \mathbb{R}$ satisfying the following conditions.
\begin{enumerate}
	\item[(i)] $\norm{\cdot} \geq 0$ \ and is $0$ only at $0$.
	\item[(ii)] $\norm{X}  =  \norm{-X}$ \ for all $X \in \mathfrak{g}$.
	\item[(iii)] $\norm{\alpha_t(X)} = \ |t|\norm{X}$ \ for all $t \in \mathbb{R}$ and $X \in \mathfrak{g}$.
	\item[(iv)] $\norm{X + Y} \leq \norm{X} + \norm{Y}$ \ for all $X, Y \in \mathfrak{g}$. 
\end{enumerate}  

Any graded Lie algebra $\mathfrak{g}$ possesses a natural candidate for a homogeneous norm as follows:
for $X = (v_1, \dots, v_r)$, let
$$\norm{X} \ = \ (\norm{v_1}_1^{2r}+\norm{v_2}_2^{2r-2}+ \cdots + \norm{v_r}_r^2)^{1/2r},$$
where $\norm{\cdot}_i$ is the Euclidean norm on each $V_i$. This norm satisfies properties (i) and (ii). However, the the homogeneity property (iii) holds if and only if $r \leq 2$. In \cite{Moskowitz1}, it was shown that the subadditivity property (iv) (or the triangle inequality) holds when $r = 2$ and the case when $r = 1$ is just the Cauchy Schwartz inequality. Our main result is that the triangle inequality holds for all $r \ge 1$.

\begin{thm}\label{mainthm}
	$\norm{X + Y} \leq \norm{X} + \norm{Y}$ \ for \ $r \geq 1$.
\end{thm}

\smallskip

\section{Proof of the main theorem}

As mentioned, the main tools we use are  H\"older's inequality and Muirhead's inequality. We first revisit here the latter inequality in the case of two variables which will be used in the proof of the main theorem.

Let $a_1, a_2, b_1$ and $b_2$ be nonnegative real numbers satisfying the following conditions:
\begin{enumerate}
	\item $a_1 \ge a_2$ and $b_1 \ge b_2$
	\item $a_1 \ge b_1$
	\item $a_1 + a_2 = b_1 + b_2$.
\end{enumerate}
Then $(a_1, a_2)$ is said to {\it majorize} $(b_1, b_2)$.

\begin{lemma}\label{Muirhead} \cite{CKP}
	Suppose that $(a_1, a_2)$ majorizes $(b_1, b_2)$. Then for nonnegative real numbers $x_1$ and $x_2$,
	\[
	x_1^{a_1}x_2^{a_2} + x_1^{a_2}x_2^{a_1} \ge x_1^{b_1}x_2^{b_2} + x_1^{b_2}x_2^{b_1}.
	\]
\end{lemma}

\vspace{0.5cm}
Now, we can prove the main theorem.
\vspace{0.5cm}

\noindent {\bf Proof of Theorem 1:}
	For convenience, we reverse the indices of summands in $X$ and $Y$. Let $X = (v_r, v_{r - 1}, \ldots, v_1)$ and $Y = (w_r, w_{r - 1}, \ldots, w_1)$. (Note that $v_r, w_r \in V_1, v_{r - 1}, w_{r - 1} \in V_2, \ldots, v_1, w_1 \in V_r$.) Then 
	\[
	X + Y = (v_r + w_r, v_{r - 1} + w_{r - 1}, \ldots,  v_1 + w_1).
	\]
	By taking $2r$th power, we need to prove that
	\[
	\norm{v_r + w_r}^{2r} + \norm{v_{r - 1} + w_{r - 1}}^{2r - 2} + \cdots + \norm{v_1 + w_1}^{2} \leq (\norm{X} + \norm{Y})^{2r}.
	\]
	Applying the Cauchy-Schwartz inequality, it is sufficient to show that 
	\[
		(\norm{v_r} + \norm{w_r})^{2r} + (\norm{v_{r - 1}} + \norm{w_{r - 1}})^{2r - 2} + \cdots + (\norm{v_1} + \norm{w_1})^2 \leq (\norm{X} + \norm{Y})^{2r}.
	\]
	Expanding both sides using the binomial theorem, the above inequality becomes
	\begin{align*}
	\sum_{k = 0}^{2r} \binom{2r}{k} \norm{v_r}^{2r - k} \norm{w_r}^k + \sum_{k = 0}^{2r - 2} \binom{2r - 2}{k}& \norm{v_{r - 1}}^{2r - 2 - k} \norm{w_{r - 1}}^k + \cdots + \sum_{k = 0}^{2} \binom{2}{k} \norm{v_1}^{2 - k} \norm{w_1}^k \\
	&\leq \\
	\sum_{k = 0}^{2r} \binom{2r}{k} &\norm{X}^{2r - k}\norm{Y}^k,   \label{eq} \tag{$\star$}
	\end{align*}
	where 
	\[
	\norm{X} = (\norm{v_r}^{2r} + \norm{v_{r - 1}}^{2r - 2} + \cdots + \norm{v_1}^2)^{1/2r} 
	\]
	and
	\[
	\norm{Y} = (\norm{w_r}^{2r} + \norm{w_{r - 1}}^{2r - 2} + \cdots + \norm{w_1}^2)^{1/2r}. 
	\]
	It is clear that the first and the last terms of each summation in (\ref{eq}) can be canceled and the inequality
	\begin{align*}
	&\binom{2r}{r - 1}\norm{v_r}^{r + 1}\norm{w_r}^{r - 1} \hspace{4cm} \binom{2r}{r - 1}\norm{X}^{r + 1}\norm{Y}^{r - 1} \\ &\hspace{2cm} + \hspace{3.5cm}  \le \hspace{4.5cm}  + \\
	&\binom{2r}{r + 1}\norm{v_r}^{r - 1}\norm{w_r}^{r + 1} \hspace{4cm} \binom{2r}{r + 1}\norm{X}^{r - 1}\norm{Y}^{r + 1}
	\end{align*} 
	holds because $\norm{v_r} \leq \norm{X}$ and $\norm{w_r} \leq \norm{Y}$.
	Now, we
	consider the middle terms of each summation in (\ref{eq}), we have
	\[
	 \binom{2r}{r}\norm{v_r}^r\norm{w_r}^r + \binom{2r - 2}{r - 1}\norm{v_{r - 1}}^{r - 1}\norm{w_{r - 1}}^{r- 1} + \cdots + \binom{2}{1}\norm{v_1}\norm{w_1} \le \binom{2r}{r}\norm{X}^r\norm{Y}^r
	\]
	which, in fact, holds by H\"older's inequality (with $p = q = 2$) and the largest coefficient $\dbinom{2r}{r}$ dominates the others. 
	
	It remains to show that for all $k = 1, 2, \ldots, r - 2$, 
	\begin{gather*}
	\left(\binom{2r}{k}\norm{v_r}^{2r - k}\norm{w_r}^k + \binom{2r}{2r - k}\norm{v_r}^{k}\norm{w_r}^{2r - k}\right) \\
	 + \left(\binom{2r - 2}{k}\norm{v_{r - 1}}^{2r - 2 - k}\norm{w_{r - 1}}^k + \binom{2r - 2}{2r - 2 - k}\norm{v_{r - 1}}^{k}\norm{w_{r - 1}}^{2r - 2 - k}\right) \\
	 \vdots\\
 + \left(\binom{2k + 2}{k}\norm{v_{k  + 1}}^{k + 2}\norm{w_{k + 1}}^k + \binom{2k + 2}{k + 2}\norm{v_{k + 1}}^{k}\norm{w_{k + 1}}^{k + 2}\right) \\
	 \le
	 \binom{2r}{k}\norm{X}^{2r - k}\norm{Y}^k + \binom{2r}{2r - k}\norm{X}^{k}\norm{Y}^{2r - k}.
	\end{gather*}
	Since the coefficient $\dbinom{2r}{k} = \dbinom{2r}{2r - k}$ is larger than the other coefficients in this inequality, we only need to show that
	\begin{align*}
		\left(\norm{v_r}^{2r - k}\norm{w_r}^k + \norm{v_r}^{k}\norm{w_r}^{2r - k}\right) 
	+ \left(\norm{v_{r - 1}}^{2r - 2 - k}\norm{w_{r - 1}}^k + \norm{v_{r - 1}}^{k}\norm{w_{r - 1}}^{2r - 2 - k}\right) \\
	+ \cdots + \left(\norm{v_{k + 1}}^{k + 2}\norm{w_{k + 1}}^k + \norm{v_{k + 1}}^{k}\norm{w_{k + 1}}^{k + 2}\right) 
	\le
	\norm{X}^{2r - k}\norm{Y}^k + \norm{X}^{k}\norm{Y}^{2r - k}.
	\end{align*}
	By H\"older's inequality (with $p = \frac{2r}{2r - k}, q = \frac{2r}{k}$ and vice versa), we have the following inequalities:
	\begin{align*}
	&\norm{v_r}^{2r - k}\norm{w_r}^k + \norm{v_{r - 1}}^{\frac{(r - 1)(2r - k)}{r}}\norm{w_{r - 1}}^{\frac{(r - 1)k}{r}} \\&\hspace{6cm} + \ldots + \norm{v_{k + 1}}^{\frac{(k + 1)(2r - k)}{r}}\norm{w_{k + 1}}^{\frac{(k + 1)k}{r}} \leq \norm{X}^{2r - k}\norm{Y}^k \ \text{and} \\
	&\norm{v_r}^{k}\norm{w_r}^{2r - k} + \norm{v_{r - 1}}^{\frac{(r - 1)k}{r}}\norm{w_{r - 1}}^{\frac{(r - 1)(2r - k)}{r}} \\&\hspace{6cm} + \ldots + \norm{v_{k + 1}}^{\frac{(k + 1)k}{r}}\norm{w_{k + 1}}^{\frac{(k + 1)(2r - k)}{r}}  \leq \norm{X}^{k}\norm{Y}^{2r - k}.
	\end{align*}
	Thus, it suffices for every $i = k + 1, k + 2, \ldots, r - 1$ to show that 
	\begin{align*}
	\norm{v_i}^{2i - k}\norm{w_{i}}^k + \norm{v_i}^{k}\norm{w_i}^{2i - k} 
	\leq
	\norm{v_i}^{\frac{i(2r - k)}{r}}\norm{w_i}^{\frac{ik}{r}} + 
	\norm{v_i}^{\frac{ik}{r}}\norm{w_i}^{\frac{i(2r - k)}{r}}. 
		\end{align*}
	For each $i = k + 1, k + 2, \ldots, r - 1$, $\left( \dfrac{i(2r - k)}{r}, \dfrac{ik}{r} \right)$ majorizes $(2i - k, k)$. By Lemma \ref{Muirhead}, we have these inequalities. \hspace*{\fill} \qedsymbol

We conclude the paper by stating the real version of Theorem \ref{mainthm} as follows.
\begin{cor}
Let $r \geq 1 $ be a positive integer, and $a_1, a_2, \ldots, a_r, b_1, b_2, \ldots, b_r$ be positive real numbers. Then 
\[
\left(\sum_{i = 1}^r (a_i + b_i)^{2i}  \right)^{\frac{1}{2r}} \leq \left(\sum_{i = 1}^r a_i^{2i}\right)^{\frac{1}{2r}} + \left(\sum_{i = 1}^r b_i^{2i}\right)^{\frac{1}{2r}}.
\]
\end{cor}
	         
\section{Acknowledgements}
The authors would like to thank the referee for suggestions that greatly improved the paper. The first author acknowledges the financial support provided by the Center of Excellence in Theoretical and Computational Science (TaCS-CoE), Faculty of Science, KMUTT.


\bigskip

\begin{thebibliography}{10}
	
	\bibitem{CKP}
	C.P. Niculescu and L.E. Persson, {\it Convex functions and their applications}, A comtemporary
	approach, CMS Books in Mathematics, vol. 23, Springer Verlag, New York, 2006.
	\bibitem{Moskowitz1}
	M. Moskowitz, \textit{An extension of Minkowski's convex body theorem to certain simply connected nilpotent groups}, Portugaliae Mathematica vol. \textbf{67}, no. 4 (2010), 541--546.
	
	\bibitem{Moskowitz2}
	M. Moskowitz, \textit{The triangle inequality for graded real vector spaces of length 3 and 4}, Math. Inequal. Appl. \textbf{17}, no. 3 (2014), 1027--1030.
	
	\bibitem{Steele}
	J. M. Steele, \textit{The Cauchy-Schwarz master class}, Cambridge University Press, 2004.
	
\end{thebibliography}
\end{document}